\documentclass[12pt]{article}

\usepackage{graphicx}%
\usepackage{multirow}%
\usepackage{amsmath,amssymb,amsfonts}%
\usepackage{amsthm}%
\usepackage{mathrsfs}%
\usepackage[title]{appendix}%
\usepackage{xcolor}%
\usepackage{textcomp}%
\usepackage{manyfoot}%
\usepackage{booktabs}%
\usepackage{algorithm}%
\usepackage{algorithmicx}%
\usepackage{algpseudocode}%
\usepackage{listings}%
\usepackage[top=2cm, bottom=2cm, right=2cm, left=2cm]{geometry}

\newtheorem{definition}{Definition}
\newtheorem{proposition}{Proposition}
\newtheorem{sled}{Corollary}

\newtheorem*{lemma}{Lemma}
\newtheorem{theor}{Theorem}
\newtheorem*{remark}{Remark}
\newcommand{\sign}{\operatorname{sign}}
\newcommand{\GL}{\operatorname{GL}}
\newcommand{\SL}{\operatorname{SL}}

\theoremstyle{definition}
\date{}
\begin{document}
\title {Essential Semigroups and Branching Rules}
\author{Andrei Gornitskii}

\maketitle

\abstract{
Let $\mathfrak{g}$ be a semisimple complex Lie algebra of finite dimension and $\mathfrak{h}$ be a semisimple subalgebra. We present an approach to find the branching rules for the pair $\mathfrak{g}\supset\mathfrak{h}$. According to an idea of Zhelobenko  the information on restriction to $\mathfrak{h}$ of all irreducible representations of $\mathfrak{g}$ is contained in one associative algebra, which we call the \emph{branching algebra}. We use an \emph{essential semigroup} $\Sigma$, which parametrizes some  bases in every finite-dimensional irreducible representations of $\mathfrak{g}$, and describe the branching rules for $\mathfrak{g}\supset\mathfrak{h}$ in terms of a certain subsemigroup $\Sigma'$ of $\Sigma$. If $\Sigma'$ is finitely generated, then the semigroup algebra corresponding to $\Sigma'$ is a toric degeneration of the branching algebra. We propose the algorithm to find a description of $\Sigma'$ in this case. We give examples by deriving the branching rules for $A_n\supset A_{n-1}$, $B_n\supset D_n$, $G_2\supset A_2$, $B_3\supset G_2$, and $F_4\supset B_4$.}

\section{Introduction}
It is an important problem in the representation theory of semisimple Lie algebras to describe how irreducible representations of a semisimple Lie algebra $\mathfrak{g}$ decompose when restricted to a semisimple Lie subalgebra $\mathfrak{h}\subset\mathfrak{g}$. This problem is called \emph{the branching problem}, and its solution is called \emph{the branching rule} or \emph{the branching law}.

If a highest weight of $\mathfrak{g}$ is fixed, then there are a number of approaches to solve branching problem for this particular highest weight. Most of them use Weyl's character formula \cite{[CKY]}, \cite{[Q]}. Another possibility is the notion of \emph{index} of the representation \cite{[P]}, \cite{[O]}. There is a computer program that derives the branching rule for a fixed highest weight \cite{[WdG]}. 

It was noticed by Zhelobenko \cite{[Z]} that if the highest weight is not fixed, then one can consider the branching problem of $\mathfrak{g}\supset\mathfrak{h}$ for all highest weights of $\mathfrak{g}$ simultaneously. This leads to the notion of \emph{branching algebra} and gives additional possibilities to study branching problems. The branching algebra collects information about branching rules for all highest weights of $\mathfrak{g}$ in one algebraic structure. This approach was used by Howe, Tan and Willenbring \cite{[HTW]} to describe the branching problems for classical symmetric pairs.
 
On the other hand, the example of Gelfand-Tsetlin bases in the irreducible representations of $\GL_n$ (or $\SL_n$) shows that the branching problem is closely related to some bases in the irreducible representations with ``good'' properties: the Gelfand-Tsetlin patterns for $SL_n$ form a finitely generated semigroup $\Sigma$ with respect to addition. This semigroup is generated by the patterns corresponding to the fundamental representations of $SL_n$. Thus the combinatorial objects parametrizing basis vectors have an additional structure. If the embedding $\mathfrak{g}\supset\mathfrak{h}$ is regular then the above connection is shown to be important in \cite{[M]}.

Our approach is a combination of these ideas and very close to the one in \cite{[M]}. Vinberg  (\emph{On some canonical bases of representation spaces of simple Lie algebras}, Conference Talk, Bielefeld, 2005) suggested a method to construct some bases parametrized by \emph{essential semigroup} $\Sigma$ in all irreducible representations of $\mathfrak{g}$. We consider the algebra $A:=\mathbb{C}[G/U]$, where $G$ is the simply connected complex algebraic group with $\mathop{\mathrm{Lie}}G=\mathfrak{g}$, and $U$ is a maximal unipotent subgroup of $G$. Then $A=\bigoplus_{\lambda} V(\lambda)^*$ is a sum of all finite-dimensional irreducible representations of $\mathfrak{g}$. We describe a subsemigroup $\Sigma'\subset\Sigma$ which parametrizes the lowest weight vectors of $\mathfrak{h}$ in $A$. If the semigroup $\Sigma'$ is finitely generated then the description of this semigroup in terms of generators and relations solves the branching problem. We give a computational algorithm to obtain such a description in this case and conjecture that $\Sigma'$ is always finitely generated. We show how this approach works, deriving the branching rules for $G_2\supset A_2$, $B_3\supset G_2$, and $F_4\supset B_4$. The last branching rule answers the question in \cite{[CKY]}. 

The semigroup $\Sigma$ is explicitly described for Lie algebras of types $A,B,C,D,G_2$ (see \cite{[FFL1]}, \cite{[FFL2]}, \cite{[G1]}, \cite{[G2]}). We prove (see Theorem \ref{theor}) that $\Sigma$ and $\Sigma'$ are closely related (under some assumptions). This connection and the description of $\Sigma$ allows one to obtain classic and well known branching rules: $A_{n}\supset A_{n-1}$, $C_n\supset C_{n-1}$, $B_n\supset D_n$. We give an example for $B_n\supset D_n$. 
\begin{remark} \emph{Vinberg's bases is a special case of \emph{essential bases} introduced by Fang,  Fourier, and Littelmann \cite{[Li]}. Our approach is suitable for this more general settings without any changes. For simplicity we will use Vinberg's approach.}
\end{remark}

 \section{Bases in the irreducible representations of $\mathfrak{g}$}
 
 All definitions and results of this section are due to Vinberg.
 \subsection{Essential semigroups}
 We recall the notion of essential semigroup that parametrizes certain bases in irreducible finite-dimensional representations of a simple complex Lie algebra.

 Let $\mathfrak{g}$ be a simple Lie algebra with the triangular decomposition
$\mathfrak{g}=\mathfrak{u}^{-}\oplus\mathfrak{t}\oplus\mathfrak{u}$, where
$\mathfrak{u}^{-}$ and $\mathfrak{u}$ are mutually opposite maximal unipotent subalgebras, and $\mathfrak{t}=\mathfrak{t}_{\mathfrak{g}}$ is a Cartan subalgebra.

One has: $\mathfrak{u}= \langle e_{\alpha}$ $\mid$ $\alpha$ $\in$ $\Delta_{+}\rangle$,
$\mathfrak{u^{-}}=\langle e_{-\alpha}$ $\mid$ $\alpha$ $\in$ $\Delta_{+}\rangle$, where
$\Delta_{+}=\Delta_+{(\mathfrak{g})}$ is the system of positive roots, $e_{\pm\alpha}$
 are the root vectors, and the symbol $\langle\ldots\rangle$ stands for the linear span.

 We denote the finite-dimensional irreducible $\mathfrak{g}$-module with highest weight $\lambda$ by $V(\lambda)$ or $V_{\mathfrak{g}}(\lambda)$ and a highest weight vector in this module by $v_{\lambda}$. We fix an ordering of positive roots : $\Delta_{+}=\{\alpha_{1},\dots,\alpha_{N}\}$.

\begin{definition}
 A signature is an $(N+1)$-tuple
 $\sigma=(\lambda;p_{{1}},\dots,p_{{N}})$,
  where $\lambda$ is a dominant weight, and $p_{i}\in\mathbb{Z}_{+}$.
\end{definition}

Set
$$
{v}(\sigma)=e_{-\alpha_{1}}^{p_{{1}}}\cdot \ldots \cdot
 e_{-\alpha_{N}}^{p_{{N}}}\cdot{v}_{\lambda}\in V(\lambda).
$$
$\lambda$ is called the \emph{highest weight} of $\sigma$, the eigenweight $\lambda-\sum p_i\alpha_i$ of the vector $v(\sigma)$ is called the \emph{weight} of $\sigma$, and the numbers $(p_1,\ldots,p_N)$ are called \emph{exponents} of $\sigma$.

 Fix any monomial order $<$ on $\mathbb{Z}^N$. We use this order to compare signatures with the same highest weight $\lambda$, i.e. if $\sigma=(\lambda;p_1,\ldots,p_N)$ and $\tau=(\lambda;q_1,\ldots,q_N)$, then $$\sigma<\tau \iff (p_1,\ldots,p_N)<(q_1,\ldots,q_N).$$

 \begin{definition} A signature $\sigma$ is essential if
 $v(\sigma)\notin\langle v(\tau)\mid\tau<\sigma\rangle$.

\end{definition}

For a dominant weight $\lambda$ the set $\{{v}(\sigma)\mid\sigma \mbox{ is essential of highest weight $\lambda$}\}$ is a basis of $V(\lambda)$. Moreover, the set of essential signatures (for all $\lambda$) is a subsemigroup of $\Lambda^+\times\mathbb{Z}_{+}^N$, where $\Lambda^+$ is the semigroup of dominant weights. The proof will be given below. We denote the semigroup of essential signatures, or \emph{essential semigroup}, by $\Sigma$.
\begin{remark}
\emph{For two weights $\lambda$ and $\mu$ we will use the notation $\lambda\succ\mu$ if $\lambda-\mu$ is a sum of positive roots. In the following we will assume that $\sigma<\tau$ provided $\lambda\succ\mu$, where $\lambda$ and $\mu$ are the weights of $\sigma$ and $\tau$, respectively. This assumption is not restrictive, because the semigroup $\Sigma$ does not depend on how we compare signatures with different weights (see the definition of essential signature). }
\end{remark}

\subsection{U-invariant functions}

Let $G$ be a simply connected simple complex algebraic group such that $\mathop{\mathrm{Lie}}G=\mathfrak{g}$.  Let $T$ be the maximal torus in $G$ such that $\mathop{\mathrm{Lie}}T=\mathfrak{t}$ and $U$ be the maximal unipotent subgroup of $G$ such that
 $\mathop{\mathrm{Lie}}U=\mathfrak{u}$.
 
Now we show that the essential signatures can be interpreted as least terms of functions on the homogeneous space $G/U$. As a consequence we prove that essential semigroup $\Sigma$ is indeed a semigroup.

 Consider the homogeneous space $G/U$.
 Let $B=T\rightthreetimes U$ be the Borel subgroup. Then

$$
 \mathbb{C}[G/U]=\bigoplus_{\lambda} \mathbb{C}[G]_{\lambda}^{(B)},
$$
where $$\mathbb{C}[G]_{\lambda}^{(B)}=\{f\in\mathbb{C}[G]\mid f(gtu)=\lambda(t)f(g),\, \forall g\in G, t\in T, u\in U\}$$
is the subspace of eigenfunctions of weight $\lambda$ for $B$ acting on $\mathbb{C}[G]$ by right translations of an argument.
 Each subspace $\mathbb{C}[G]_{\lambda}^{(B)}$ is finite-dimensional and is isomorphic as a $G$-module (with respect to the action of $G$ by left translations of an argument), to the space $V(\lambda)^{*}$ of linear functions on $V(\lambda)$ (see \cite{[Ï]}, Theorem 3). The isomorphism is given by the formula:
$$
 V(\lambda)^{*}\ni\omega \longmapsto f_{\omega}\in\mathbb{C}[G]_{\lambda}^{(B)},
 \quad\textrm{where}\quad f_{\omega}(g)=\langle\omega, g\emph{v}_{\lambda}\rangle.
$$

Let $U^{-}$ be the maximal unipotent subgroup such that $\mathop{\mathrm{Lie}}U^{-}=\mathfrak{u^{-}}$. The function $f_{\omega}$ is uniquely determined by its restriction to the dense open subset $U^{-}\cdot$$T\cdot$$U$; moreover
\begin{multline*}
 $$f_{\omega}(u^{-}\cdot t\cdot u)=\langle\omega
,u^{-} tu\emph{v}_{\lambda}\rangle=\langle\omega
,\lambda(t)u^{-}\emph{v}_{\lambda}\rangle=\lambda(t)f_{\omega}(u^{-}),\\ \quad \forall u\in U, u^{-}\in U^{-}, t\in T.
$$
\end{multline*}
Next, $U^{-}=U_{-\alpha_{1}}\cdot\ldots\cdot U_{-\alpha_{N}}$, where
$U_{\alpha}=\{\exp(ze_{\alpha})\mid$ $z\in\mathbb{C}\}$ (see \cite[Sec. X, \S 28.1]{[X]}). Hence

$$
 u^{-}=\exp(z_{1}e_{-\alpha_{1}})\cdot\ldots\cdot\exp(z_{N}e_{-\alpha_{N}}).
$$
Thus we obtain
$$
f_{\omega}(u^{-})=\left\langle\omega,
\exp(z_{1}e_{-\alpha_{1}})\cdot\ldots\cdot\exp(z_{N}e_{-\alpha_{N}})
\cdot \emph{v}_{\lambda}\right\rangle=\sum_{\sigma=(\lambda;p_{1},\dots,p_{N})}
\frac{\prod z_{i}^{p_{i}}}{\prod p_{i}!}
\langle\omega,\emph{v}(\sigma)\rangle.
$$

\begin{proposition}
 A signature $\sigma$ is essential if and only if
$\prod z_{i}^{p_{i}}$ is the least term in
 $f_{\omega}|_{U^{-}}$ for some
$\omega \in V(\lambda)^{*}$ in the sense of the order introduced above.
\end{proposition}
\begin{proof}
Let $\prod z_{i}^{p_{i}}$ be the least term in $f_{\omega}|_{U^{-}}$ for some $\omega\in V(\lambda)^{*}$. Then $\omega$ vanishes on all vectors $\emph{v}(\tau)$ with $\tau<\sigma$ and is nonzero at  $\emph{v}(\sigma)$. Consequently, $\emph{v}(\sigma)$ cannot be expressed via $\emph{v}(\tau)$ with $\tau<\sigma$, and hence  $\sigma$ is essential.

Conversely, let $\sigma$ be essential. Consider a function $\omega$ that vanishes on $\emph{v}(\tau)$ for all essential $\tau$ except for $\sigma$. Obviously, $f_{\omega}|_{U^{-}}$ has the desired least term.
\end{proof}
\begin{proposition}
If $\sigma,\tau\in\Sigma$ then $\sigma+\tau\in\Sigma$.
\end{proposition}
\begin{proof}
Suppose that the least terms in $f|_{U^{-}}$ and $g|_{U^{-}}$ correspond to the essential signatures $\sigma$ and $\tau$. Then the least term in $(f\cdot g)|_{U^{-}}$ corresponds to the signature $\sigma+\tau$. Hence $\sigma+\tau$ is essential.
\end{proof}

\section{The branching algebra and the branching semigroup}
\subsection{The branching problem and the branching algebra}
Let $\mathfrak{h}\subset\mathfrak{g}$ be a simple Lie subalgebra of $\mathfrak{g}$. Let $H\subset G$ be a connected algebraic group such that $\mathop{\mathrm{Lie}}H=\mathfrak{h}$.

 Restrict the irreducible representation $V_{\mathfrak{g}}(\lambda)$ with the highest weight $\lambda$ to $\mathfrak{h}$:
$$
V_{\mathfrak{g}}(\lambda)|_{\mathfrak{h}}=\bigoplus_{\lambda'}m_{\lambda,\lambda'}V_{\mathfrak{h}}(\lambda'),
$$
where $V_{\mathfrak{h}}(\lambda')$ is the irreducible representation of $\mathfrak{h}$ with the highest weight $\lambda'$, and $m_{\lambda,\lambda'}$ is the multiplicity. The classical branching problem is to determine $m_{\lambda,\lambda'}$.

Consider the action of $H$ on $\mathbb{C}[G/U]=\bigoplus_{\lambda}V(\lambda)^*$ by left translations of an argument. Let $U'$ be a maximal unipotent subgroup of $H$ such that $\mathop{\mathrm{Lie}}U'=\mathfrak{u'}$. The algebra $\mathbb{C}[G/U]^{U'}$ of $U'$-invariants is called \emph{the branching algebra}. This is a finitely generated algebra consisting of the highest vectors of $\mathfrak{h}$. A description of this algebra in terms of generators and relations solves the branching problem.

\subsection{The branching semigroup}

Now we want to introduce a subsemigroup $\Sigma'$ of the essential semigroup $\Sigma$, which is related to the branching problem. We call $\Sigma'$ \emph{the branching semigroup}.

Recall that $f_\omega\in V(\lambda)^*$ is uniquely determined by its restriction to $U^-\cdot T$. Let $t_1,\ldots,t_n$ be the coordinates on $T$ corresponding to the fundamental weights $\pi_i$, i.e. $t_i=\pi_i(t), t\in T$. Then $f_\omega$ can be thought as a polynomial in $t_1,\ldots,t_n, z_1,\ldots, z_N$. Indeed, if $\lambda=\sum_i k_i\pi_i$ then $$f_\omega(u^-\cdot t)=t_1^{k_1}\cdot\ldots\cdot t_n^{k_n}\cdot \left(\sum_{\sigma=(\lambda;p_{1},\dots,p_{N})}
\frac{\prod z_{i}^{p_{i}}}{\prod p_{i}!}
\langle\omega,\emph{v}(\sigma)\rangle\right).$$
The expression in the brackets has the form $cz_1^{p_1}\cdot\ldots\cdot z_N^{p_N}+\mbox{higher terms}$, where $c\in\mathbb{C}\backslash \{0\}$. Set $\sign(f_\omega)=(\lambda;p_1,\ldots,p_N)\in\Sigma$. Obviously, $\sign(f_{\omega_1} f_{\omega_2})=\sign(f_{\omega_1})+\sign(f_{\omega_2})$.

Let $\Sigma'=\{\sign(f_\omega)\mid f_\omega\in\mathbb{C}[G/U]^{U'^{-}}\}$, where $U'^{-}\subset H$ is the opposite maximal unipotent subgroup to $U'$.  So, $\Sigma'$ consists of essential signatures that are the least terms of the lowest vectors with respect to $\mathfrak{h}$. Denote by $\Sigma'(\lambda)$ the set of all signatures of the highest weight $\lambda$ in $\Sigma'$.

 If $\Sigma'$ is finitely generated then a description of $\Sigma'$ in terms of generators and relations solves the branching problem. Indeed, the signature $\sigma\in\Sigma'(\lambda)$ defines the irreducible representation $V_{\mathfrak{h}}(\lambda')$ in $V_{\mathfrak{g}}(\lambda)$ where $\lambda'$ is the weight of $v(\sigma)$ restricted to $\mathfrak{h}$. Therefore the multiplicity
  $m_{\lambda,\lambda'}$ is equal to the number of signatures $\sigma$ in $\Sigma'(\lambda)$ such that the weight of $v(\sigma)$ is $\lambda'$ when restricted to $\mathfrak{h}$.

\subsection{ Approaches to solve the branching problem}
In this section we discuss computational approach and theoretical approach to describe the semigroup $\Sigma'$. 

The computational approach is straightforward and works if and only if $\Sigma'$ is finitely generated. We conjecture that this is always the case. Despite the simplicity, this method is usefull. In the next section we give examples deriving the branching rules for $G_2\supset A_2$ and $B_3\supset G_2$ in probably the simplest known way. Moreover, this method allows to obtain the branching rule for $F_4\supset B_4$, that answers the question in \cite{[CKY]}. Finally, we derive the branching rule for $A_n\supset A_{n-1}$ by using the combinatorial result on the number of semistandard Young tableaux.
 
 We will need the following lemma:
\begin{lemma} Let $f,g\in\mathbb{C}[x_0,\ldots,x_n]$ be polynomials of total degree $k$, and let $f(\lambda)=g(\lambda)$ for all $\lambda\in\{(\lambda_0,\ldots,\lambda_n)\in\mathbb{Z}_+^{n+1}\mid \lambda_0+\ldots+\lambda_n\leq k\}:=I.$ Then $f=g$. 
\end{lemma}
\begin{proof}
Let $f=\sum_{\lambda\in I}f_\lambda x^\lambda$ and $g=\sum_{\lambda\in I}g_\lambda x^\lambda$, where $f_\lambda,g_\lambda\in\mathbb{C}, x^{\lambda}=x_0^{\lambda_0}\cdot\ldots\cdot x_n^{\lambda_n}$. We know that $\sum_{\lambda\in I} (f_\lambda- g_\lambda)\mu_0^{\lambda_0}\cdot\ldots\cdot\mu_n^{\lambda_n}=0$ for all $\mu\in I$. So we obtain a system of linear equations with respect to the variables $f_\lambda-g_\lambda$. To prove that $f=g$ it is enough to prove that the corresponding $\mid I\mid\times\mid I\mid$-matrix $(\mu^\lambda), \mu,\lambda\in I$, is invertible. This follows directly from Theorem 1 in \cite{[BCR]} after substitution $x_{i,j}=j$.
\end{proof}
The computational approach is based on the following steps:
\begin{enumerate}
\item [step 1:] Choose some set of dominant weights $S=\{\lambda_1,\ldots,\lambda_s\}$ including the fundamental weights. Find the sets $\Sigma'(\lambda_i)$ and generate a semigroup $\Sigma'_{S}$ by $\Sigma'(\lambda_i)$. \\
\item [step 2:]  Compute the sum $$d(\lambda):=\sum_{\sigma=(\lambda;\ldots)\in\Sigma'_{S}} \dim V_{\mathfrak{h}}(\lambda').$$ The sum is taken over the signatures in $\Sigma'_S$ with the highest weight $\lambda$. $V_{\mathfrak{h}}(\lambda')$ is the irreducible representation of $\mathfrak{h}$ corresponding to the signature $\sigma\in\Sigma'_{S}$.\\
\item [step 3:] If $d(\lambda)=\dim V_{\mathfrak{g}}(\lambda)$ for sufficiently large number of $\lambda=\sum t_j\pi_j$ (see the above lemma), and $d(\lambda)$ is a polynomial in $t_j$ then $\Sigma'(\lambda_i)$ generate $\Sigma'$. 
Otherwise, the equality $d(\lambda)=\dim V_{\mathfrak{g}}(\lambda)$ fails for some $\lambda$. Add $\lambda$ to $S$ and repeat the steps 1-3.  
\end{enumerate}

In step 3 we use the fact that $\dim V_{\mathfrak{g}}(\lambda)$, where $\lambda=\sum t_j\pi_j$, is a polynomial in $t_j$ by Weyl's dimension formula. Obviously, the approach above is an algorithm if and only if the semigroup $\Sigma'$ is finitely generated.

The theoretical approach of describing $\Sigma'$ is based on the connection with $\Sigma$. The semigroup $\Sigma$ is explicitly described for Lie algebras of types $A,B,C,D,G_2$ (see \cite{[FFL1]}, \cite{[FFL2]}, \cite{[G1]}, \cite{[G2]}) for some ordering of positive roots and monomial order. 

Let $\mathfrak{h}\subset\mathfrak{g}$ be a regular embedding such that $\mathfrak{t}_\mathfrak{h}\subset \mathfrak{t}_\mathfrak{g}$ and $\Delta_+(\mathfrak{h})\subset\Delta_+(\mathfrak{g})$. Let $\tilde{\Delta}_+:=\Delta_+(\mathfrak{g})\backslash\Delta_+(\mathfrak{h})$. For a signature $\sigma=(\lambda;p_1,\ldots,p_N)$ denote by $\tilde\sigma$ the signature $(\lambda;\tilde{p}_1,\dots,\tilde{p}_N)$ such that $\tilde{p}_i=p_i$ if $\alpha_i\in\tilde{\Delta}_+$ and $\tilde{p}_i=0$ if $\alpha_i\in\Delta_+(\mathfrak{h})$.

 We say that $\Sigma$ (or rather the chosen ordering of positive roots of $\mathfrak{g}$ and monomial order) is \emph{compatible} with the embedding $\mathfrak{h}\subset\mathfrak{g}$ if the following hold:
\begin{itemize}
\item[(i)] The roots $\Delta_+(\mathfrak{h})$ precede the roots $\tilde\Delta_+$ in the ordering of positive roots of $\mathfrak{g}$,
\item[(ii)] if $\tilde\sigma <\tilde\mu$ then $\sigma<\mu$. 
\end{itemize}

\begin{theor}
\label{theor}
Let the semigroup $\Sigma$ is compatible with the embedding $\mathfrak{h}\subset\mathfrak{g}$. Then $\Sigma'=\{\sigma\in\Sigma\mid\sigma=\tilde\sigma\}$.
\end{theor}

\begin{proof}
Let $\sigma=(\lambda;p_1,\ldots,p_N)\in\Sigma'$. This means that $\sigma=\sign(v_{\mu}^*)$, where $V^*_{\mathfrak{h}}(\mu)\subset V^*_{\mathfrak{g}}(\lambda)$ and $v_{\mu}^*$ is the lowest vector in $V^*_{\mathfrak{h}}(\mu)$. The signature $\sigma$  is the minimal signature satisfying $\langle v^*_{\mu},v(\sigma)\rangle\neq 0.$ The vector $v(\sigma)$ has nonzero projection $c\cdot v_{\mu}$ on $V_{\mathfrak{h}}(\mu)$, where $v_{\mu}$ is the highest vector and $c\in\mathbb{C}$. We want to show that $\sigma=\tilde\sigma$. Suppose $\sigma\neq\tilde\sigma$. Then (see (i)) $v(\sigma)=e_{-\alpha_{1}}^{p_1}\cdot\ldots\cdot e_{-\alpha_{s}}^{p_s}\cdot v(\tilde\sigma)$, where $\sum p_i > 0$ and $\alpha_1,\ldots,\alpha_s$ are the roots of $\mathfrak{h}$. Obviously, $v(\tilde\sigma)$ has zero projection on $V_{\mathfrak{h}}(\mu)$. Since the projection is $\mathfrak{h}$-invariant then $v(\sigma)$ has zero projection on $V_{\mathfrak{h}}(\mu)$ as well. A contradiction.

Conversely, let $\omega\in V^*_{\mathfrak{h}}(\mu)$ be any weight vector that is not the lowest vector. Let $\sigma=\sign(v_{\mu}^*)$, where  $v_{\mu}^*$ is the lowest vector, and
let $\sign(\omega)=\sigma_{\omega}.$

One has $\langle v(\sigma),v^*_{\mu}\rangle\neq 0$. Therefore there exists a signature $\sigma'=(\lambda;p_1,\ldots,p_s,\ldots)$ satisfying $\tilde\sigma'=\sigma$ and $v(\sigma')= e_{-\alpha_{1}}^{p_1}\cdot\ldots\cdot e_{-\alpha_{s}}^{p_s}\cdot v(\sigma)$, where $\alpha_1,\ldots,\alpha_s$ are the roots of $\mathfrak{h}$ and $\sum p_i>0,$ such that $\langle v(\sigma'),\omega\rangle\neq 0$.

 We claim that $\tilde\sigma_{\omega}=\sigma$. Indeed, the vector $v(\tilde\sigma_{\omega})$ has nonzero projection on $V_{\mathfrak{h}}(\mu)$. Hence if 
$$\textrm{the weight of }\,\tilde\sigma_{\omega} \prec \textrm{the weight of }\,\sigma,
$$
then (ii) implies that $\sigma' < \sigma_{\omega}$. A contradiction, because $\langle v(\sigma'),\omega\rangle\neq 0$ and the signature $\sigma_\omega$ is a minimal signature satisfying $\langle v(\sigma_\omega),\omega\rangle\neq 0$.
  Hence both $\tilde{\sigma}_{\omega}$ and $\sigma$ have the same weight and satisfy $\langle v(\sigma),v^*_{\mu}\rangle\neq 0$ and $\langle v(\tilde\sigma_{\omega}),v^*_{\mu}\rangle\neq 0$, respectively. Since $\sigma_\omega$ is $\sign(\omega)$ and $\sigma$ is the minimal signature satisfying $\langle v(\sigma),v^*_{\mu}\rangle\neq 0$ then (ii) implies that $\tilde\sigma_{\omega}=\sigma$. We conclude that $\sigma_{\omega}\neq \tilde\sigma_{\omega}$, because $\omega$ is not the lowest vector. So, $\sigma_{\omega}\notin\Sigma'$.
  
  Finally, for every $\mu$ we represent the isotypic component $V^*_{\mathfrak{h}}(\mu)\oplus\ldots\oplus V^*_{\mathfrak{h}}(\mu)$ of the highest weight $\mu$ in such a way that the signatures of the lowest vectors are different. It was proved above that if the signature of the lowest vector in $V^*_{\mathfrak{h}}(\mu)$ is $\sigma$, then $\widetilde{\sign(\omega)}=\sigma$ for every $\omega\in V^*_{\mathfrak{h}}(\mu)$. Hence the signatures $\sign(\omega)$ are different for different $\mu$ and different summands in the decomposition of the isotypic component of the highest weight $\mu$. Therefore for any $v^*\in V^*(\lambda)$ the signature $\sigma=\sign(v^*)$ coincide with $\sign(\omega),\omega\in V^*_{\mathfrak{h}}(\mu)$, for some $\mu$, $\omega$ and some summand in the decomposition of the isotypic component of the highest weight $\mu$. Hence if $\sigma=\tilde\sigma$ then $\sigma=\sign(v^*_{\mu})$, where $v^*_{\mu}$ is the lowest vector of some $V_{\mathfrak{h}}(\mu)$. This implies $\sigma\in\Sigma'$.
\end{proof}
\begin{sled}
\label{cor} If $\Sigma$ is generated by signatures $S:=\{\sigma_1,\ldots,\sigma_m\}$, then $\Sigma'$ is generated by signatures $\{\sigma\in S\mid \tilde\sigma=\sigma\}$.
\end{sled}
\begin{sled} If $\sigma=\tilde\sigma\in\Sigma$, then $\sigma$ is a signature corresponding to some lowest vector with respect to $\mathfrak{h}$.
\end{sled}
The description of $\Sigma$ for Lie algebras of types $A,B,C,D$ and the above theorem can be used to describe classic branching rules: $A_{n-1}\subset A_n$, $D_n\subset B_n$, $C_{n-1}\subset C_n$. In the next section we give an example for $D_n\subset B_n$.

 \section{Examples}

In this section we give examples of descriptions of the branching semigroup $\Sigma'$ for the embeddings $B_n\supset D_n$, $A_{n}\supset A_{n-1}$, $G_2\supset A_2$, $B_3\supset G_2$, and $F_4\supset B_4$. As was  noted above the description of $\Sigma'$ solves the branching problem. 

In the case $B_n\supset D_n$ we use the previous theorem and the description of $\Sigma$ given in \cite{[G2]}. In the rest cases we use the computational approach. We omit some computational details (that can be easily done by using a computer). For example we omit the computing of $d(\lambda)$ and the comparison with $\dim V(\lambda)$.
\subsection{The branching rule for $B_n\supset D_n$}
We introduce some notation and recall basic facts about representations of orthogonal Lie algebras. The numeration of fundamental weights is according to \cite[Table 1]{[VO]}.

 We denote the fundamental weights for $B_n$ by the same letters as for $D_n$, by abuse of notation. Let $\widehat\omega_p=\omega_p$ if $p\neq n-1$ and $\widehat\omega_{n-1}=\omega_{n-1}+\omega_{n}$ for $D_n$, and let
$\widehat\omega_p=\omega_p$ if $p\neq n$ and $\widehat\omega_{n}=2\omega_{n}$ for $B_n$.

Recall that $V(\omega_1)$ is the standard representation of $\mathfrak{so}_{2n+1}$ (resp. $\mathfrak{so}_{2n}$) in $\mathbb{C}^{2n+1}$ (resp. $\mathbb{C}^{2n}$).

Let $\pm\varepsilon_i$ ($i=1,\ldots, n$) be the nonzero weights of the representation $V(\omega_1)$ of $D_n$ or $B_n$. Then the positive roots of $D_n$ are
$$
\varepsilon_i\pm\varepsilon_j, \quad i<j,\quad i,j\in\{1,\ldots,n\},
$$
and the positive roots of $B_n$ are
$$
\varepsilon_i\pm\varepsilon_j, \quad i<j, \quad i,j\in\{1,\ldots,n\},
$$
$$
\varepsilon_i,\quad i\in\{1,\ldots,n\}.
$$
The fundamental weights and weights $\widehat\omega_i$ can be expressed via $\varepsilon_i$ as follows:
$$
\widehat\omega_i=\varepsilon_1+\ldots+\varepsilon_i,\quad i=1,\ldots,n \quad\textrm{for}\quad B_n,\quad i=1,\ldots,n-1\quad\textrm{for}\quad D_n;
$$

$$
\omega_n=\frac{1}{2}(\varepsilon_1+\ldots+\varepsilon_n)\quad\textrm{for both $B_n$ and $D_n$};
$$
$$
\omega_{n-1}=\frac{1}{2}(\varepsilon_1+\ldots+\varepsilon_{n-1}-\varepsilon_n)\quad\textrm{for}\quad D_n.
$$

Denote by $e_{\pm i}$ eigenvectors in $V(\omega_1)$ of eigenvalues $\pm\varepsilon_i$, and denote by $e_0$ an eigenvector of eigenvalue 0 (for $B_n$).

We have the standard embedding of $SO_{2n}$ in $SO_{2n+1}$ such that the following $D_n$-module decomposition holds:
$$
V_{B_n}(\omega_1)=V_{D_n}(\omega_1)\oplus\langle e_0\rangle.
$$

Now we recall the description of $\Sigma$ according to \cite{[G2]}.

We choose a numeration on the sets of positive roots for $B_n$ as follows:
$$
\varepsilon_1-\varepsilon_2,\varepsilon_1+\varepsilon_2,\ldots,\varepsilon_{1}-\varepsilon_n,\ldots,\varepsilon_{n-1}-\varepsilon_{n},\varepsilon_1+\varepsilon_n,\ldots,\varepsilon_{n-1}+\varepsilon_n, \varepsilon_1,\ldots,\varepsilon_n.
$$
Obviously, this numeration satisfies (i).

Also we have a monomial order on the set of signatures. We compare two signatures of $B_n$ of the same highest weight as follows (we move on to the next step if on the previous steps the tuples of exponents of the signatures coincide):
\begin{enumerate}
\item compare the tuples of exponents corresponding to the roots
 $\varepsilon_1,\ldots,\varepsilon_n$ by using the degree lexicographic order,
\item compare the tuples of exponents corresponding to the roots
 $\varepsilon_1+\varepsilon_n,\ldots,\varepsilon_{n-1}+\varepsilon_n$ by using the degree lexicographic order,

 \item compare the tuples of exponents corresponding to the roots $\varepsilon_1-\varepsilon_n,\ldots,\varepsilon_{n-1}-\varepsilon_n$ by the degree lexicographic order,
     \item compare the tuples of exponents corresponding to the roots $\varepsilon_1+\varepsilon_{n-1},\ldots,\varepsilon_{n-2}+\varepsilon_{n-1}$ by the degree lexicographic order,
         \item compare the tuples of exponents corresponding to the roots $\varepsilon_1-\varepsilon_{n-1},\ldots,\varepsilon_{n-2}-\varepsilon_{n-1}$ by the degree lexicographic order,
 \item[]\ldots\ldots
 \item compare the exponents corresponding to the root $\varepsilon_1+\varepsilon_2$,
     \item compare the exponents corresponding to the root $\varepsilon_1-\varepsilon_2$.
\end{enumerate}
This monomial order satisfies (ii). Therefore the corresponding $\Sigma$ is compatible with the embedding $D_n\subset B_n$.
The semigroup $\Sigma$ is generated by essential signatures of highest weights in the set $\{\omega_1,\ldots,\omega_n,2\omega_n\}$ (see Theorem 2 in \cite{[G2]}). One has the following decompositions:
\begin{center}
$
V_{B_n}(\omega_k)^*\mid_{D_{n}}= V_{D_{n}}(\hat{\omega}_k)^*\oplus V_{D_{n}}(\omega_{k-1})^*, k=1,\ldots,n,$\\

$V_{B_n}(2\omega_n)^*\mid_{D_{n}}= V_{D_{n}}(\hat{\omega}_{n-1})^*\oplus V_{D_{n}}(2\omega_{n-1})^*\oplus V_{D_{n}}(2\omega_{n})^*$,\\
\end{center}
where $\omega_0=0$.

For every representation $V_{D_{n}}(\mu)^*$ in $V_{B_n}(\lambda)^*$, where $\lambda=\omega_1,\ldots,\omega_n,2\omega_n$, we attach a signature $\sigma_{\lambda;\mu}$ of the lowest vector in $V_{D_{n}}(\mu)^*$. In what follows we denote by $(\lambda;0)$ the signature with the highest weight $\lambda$ and zero exponents corresponding to all positive roots, and we denote by $(\lambda;k\alpha_i)$ the signature with the highest weight $\lambda$ and with the only nonzero exponent $p_i=k$ corresponding to the positive root $\alpha_i$. One has 

\begin{flushleft}
\begin{tabular}{ll}
$\sigma_{2k-1}:=(\omega_k;0)=\sigma_{\omega_k,\hat{\omega}_k}$& $\sigma_{2k}:=(\omega_k;\varepsilon_k)=\sigma_{\omega_k,\omega_{k-1}}$, $k=1,\ldots,n,$\\
$\sigma_{2n+1}:=(2\omega_{n};\varepsilon_n)=\sigma_{2\omega_n,\hat{\omega}_{n-1}}$&$\sigma_{2n+2}:=(2\omega_n;2\varepsilon_n)=\sigma_{2\omega_n,2\omega_{n-1}}$\\
$\sigma_{2n+3}:=(2\omega_n,0)=\sigma_{2\omega_n,2\omega_n}$&
\end{tabular}

\end{flushleft}

\begin{theor} $\Sigma'$ is a free semigroup generated by $\sigma_i, i=1,\ldots, 2n$.
\end{theor}
\begin{proof}
It was already noted that the semigroup $\Sigma$ is generated by essential signatures of highest weights in the set $\{\omega_1,\ldots,\omega_n,2\omega_n\}$. Then Corollary \ref{cor} implies that $\Sigma'$ is generated by essential signatures $\sigma$ of highest weights in $\omega_1,\ldots,\omega_n,2\omega_n$ satisfying $\tilde\sigma=\sigma$. Obviously, the signatures $\sigma_{2n+1}$, $\sigma_{2n+2}$, and $\sigma_{2n+3}$ belong to the semigroup generated by $\sigma_i, i=1,\ldots, 2n$. 
\end{proof}
\subsection{The branching rule for $A_n\supset A_{n-1}$}
We denote by $\varepsilon_i, i=1,\ldots,n+1,$ the weights of the standard representation of $A_n$ in $\mathbb{C}^{n+1}$. Here $\varepsilon_i$ is a weight of the vector $e_i$ with respect to Cartan subalgebra of $A_n$ consisting of diagonal matrices. Let $\beta_i=\varepsilon_i-\varepsilon_{i+1}, i=1,\ldots,n$, be the simple roots and let $\pi_i=\varepsilon_1+\ldots+\varepsilon_i$ be the fundamental weights. Choose any ordering of positive roots of $A_n$. Choose any homogeneous order on $\mathbb{Z}^{\frac{(n+1)n}{2}}$.

Let $A_{n-1}$ consists of matrices in $A_n$ with a zero last row and column. We will denote the fundamental weights of $A_{n-1}$ by $\omega_1,\ldots,\omega_{n-1}$. 

Restrict the irreducible representations of $A_n$ to $A_{n-1}$. One has
\begin{center}
$
V_{A_n}(\pi_k)^*\mid_{A_{n-1}}= V_{A_{n-1}}(\omega_k)^*\oplus V_{A_{n-1}}(\omega_{k-1})^*, k=1,\ldots,n,
$\\
\end{center}
where $\omega_n=\omega_0=0$.

For every representation $V_{A_{n-1}}(\lambda)^*$ in $V_{A_n}(\pi_k)^*, k=1,\ldots,n,$ we attach a signature $\sigma_{\pi_k;\lambda}$ of the lowest vector in $V_{A_{n-1}}(\lambda)^*$. In what follows we denote by $(\lambda;0)$ the signature with the highest weight $\lambda$ and zero exponents corresponding to all positive roots, and we denote by $(\lambda;\alpha_i)$ the signature with the highest weight $\lambda$ and with the only nonzero exponent $p_i=1$ corresponding to the positive root $\alpha_i$. One has 

\begin{flushleft}
\begin{tabular}{ll}
$\sigma_1:=(\pi_1;0)=\sigma_{\pi_1,\omega_1}$& $\sigma_2:=(\pi_1;\varepsilon_1-\varepsilon_{n+1})=\sigma_{\pi_1,\omega_0}$\\
$\sigma_3:=(\pi_2;0)=\sigma_{\pi_2,\omega_2}$&$\sigma_4:=(\pi_2;\varepsilon_2-\varepsilon_{n+1})=\sigma_{\pi_2,\omega_1}$\\
\dotfill & \dotfill\\
$\sigma_{2n-1}:=(\pi_n;0)=\sigma_{\pi_n,\omega_n}$& $\sigma_{2n}:=(\pi_n;\varepsilon_{n}-\varepsilon_{n+1})=\sigma_{\pi_n,\omega_{n-1}}$
\end{tabular}

\end{flushleft}

\begin{theor} $\Sigma'$ is a free semigroup generated by $\sigma_i, i=1,\ldots, 2n$.
\end{theor}
\begin{proof}
Denote by $\Sigma'_S$ the semigroup generated by $\sigma_i, i=1,\ldots, 2n.$ Fix a dominant weight $\lambda=\sum_{i=1}^n k_i\varepsilon_i=\sum_{i=1}^n (k_i-k_{i+1})\pi_i$, where $k_1\geq\ldots\geq k_{n}\geq k_{n+1}=0$. It is enough to show that $\dim V_{A_n}(\lambda)=\sum_{\sigma=(\lambda;\ldots)\in\Sigma'_{S}} \dim V_{A_{n-1}}(\lambda')$, where $\lambda'$ is the highest weight with respect to $A_{n-1}$ corresponding to $\sigma$. The signatures 
$$
\sum_{i=1}^n(k_i-k_{i+1}-k'_i)\sigma_{2i-1}+k'_i\sigma_{2i},\quad 0\leq k_i'\leq k_i-k_{i+1},
$$
are all signatures in $\Sigma'_S$ with the highest weight $\lambda$. The weight of such signature is $\sum_{i=1}^{n}(k_i-k'_i)\varepsilon_i$. So we have all weights $\mu_1\varepsilon_1+\ldots+\mu_n\varepsilon_n$, where $k_1\geq\mu_1\geq k_2\geq\mu_2\geq\ldots k_n\geq\mu_n\geq 0.$ In terms of partitions these weights correspond to Young diagrams $\lambda'$ obtained from the diagram $\lambda$ by deleting at most one box in each column.

It is known that $\dim V_{A_n}(\lambda)$ is the number of semistandard Young tableaux corresponding to the partition $k_1\geq k_2\geq\ldots\geq k_n\geq 0$ filled with numbers from $1$ to $n+1$. Deleting the boxes with $n+1$ entry gives a bijection between all semistandard Young tableaux of the form $\lambda$ and all semistandard Young tableaux of the form $\lambda'$ (for all $\lambda'$) filled with numbers from $1$ to $n$. So the equality $d(\lambda)=\dim V_{A_n}(\lambda)$ holds for all $\lambda$. Hence $\Sigma'_S=\Sigma'$.
\end{proof}
  \subsection{The branching rule for $G_2\supset A_2$}

 Choose an ordering of positive roots of $G_2$ as follows:

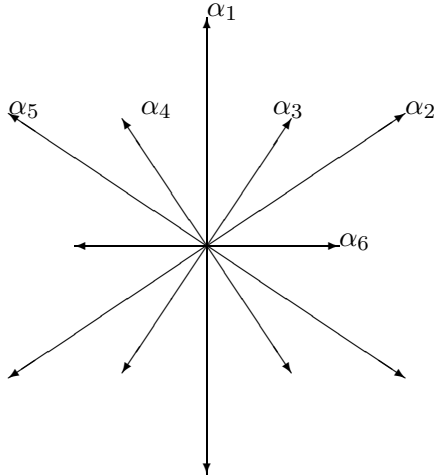
\begin{picture}(200,200)
\put(100,100){\vector(3,2){75}}
\put(175,150){$\alpha_{2}$}
\put(100,100){\vector(3,-2){75}}
\put(100,100){\vector(-3,-2){75}}
\put(100,100){\vector(0,-1){86.5}}
\put(100,100){\vector(-3,2){75}}
\put(25,150){$\alpha_{5}$}
\put(100,100){\vector(0,1){86.5}}
\put(100,187){$\alpha_{1}$}
\put(100,100){\vector(1,0){50}}
\put(150,100){$\alpha_{6}$}
\put(100,100){\vector(-1,0){50}}
\put(100,100){\vector(2,3){32}}
\put(125,150){$\alpha_{3}$}
\put(100,100){\vector(-2,3){32}}
\put(75,150){$\alpha_{4}$}
\put(100,100){\vector(2,-3){32}}
\put(100,100){\vector(-2,-3){32}}

\end{picture}

The roots of $A_2$ are the long roots of $G_2$. Choose any homogeneous order on $\mathbb{Z}^6$. This gives us the essential semigroup $\Sigma$. We will denote the fundamental weights of $G_2$ by $\pi_1,\pi_2$, and the fundamental weights of $A_2$ by $\omega_1, \omega_2$ to avoid any confusions. 

Restrict the representation of $G_2$ with the highest weight $\pi_1$ ($\dim=7$) and $\pi_2$ ($\dim=14$) to $A_2$:
\begin{center}
$
V_{G_2}(\pi_1)^*\mid_{A_2}= V_{A_2}(\omega_1)^*\oplus V_{A_2}(\omega_2)^*\oplus V_{A_2}(0)^*.
$\\

$
V_{G_2}(\pi_2)^*\mid_{A_2}=V_{A_2}(\omega_1)^*\oplus V_{A_2}(\omega_2)^*\oplus V_{A_2}(\omega_1+\omega_2)^*.
$
\end{center}
For every representation $V_{A_2}(\lambda)^*$ in $V_{G_2}(\pi_i)^*, i=1,2,$ we attach a signature $\sigma_{\pi_i;\lambda}$ of the lowest vector in $V_{A_2}(\lambda)^*$. One has 

\begin{tabular}{ll}
$\sigma_1:=(\pi_1;0,0,0,0,0,1)=\sigma_{\pi_1,\omega_1}$& $\sigma_4:=(\pi_2;0,0,0,0,0,0)=\sigma_{\pi_2,\omega_1+\omega_2}$\\
$\sigma_2:=(\pi_1;0,0,0,0,0,0)=\sigma_{\pi_1,\omega_2}$&$\sigma_5:=(\pi_2;0,0,1,0,0,0)=\sigma_{\pi_2,\omega_1}$\\
$\sigma_3:=(\pi_1;0,0,1,0,0,0)=\sigma_{\pi_1,0}$&$\sigma_6:=(\pi_2;0,0,0,1,0,0)=\sigma_{\pi_2,\omega_2}$
\end{tabular}

We check that $d(\lambda)$ is a polynomial and $d(\lambda)=\dim V(\lambda)$. So, we conclude that $\sigma_i,i=1,\ldots,6,$ generate $\Sigma'$ with one relation $\sigma_2+\sigma_5=\sigma_3+\sigma_4$. This gives the solution of the branching problem. 

\subsection{The branching rule for $B_3\supset G_2$}

The representation $V(\omega_1)$ ($\dim=7$) of $G_2$ admits a nondegenerate symmetric $G_2-$invariant bilinear form. So one has the embedding $\rho: G_2\to\mathfrak{so}_7=B_3$.

Let $\beta_1, \beta_2, \beta_3$ be the  simple roots for $B_3$ and let $\pi_1, \pi_2, \pi_3$ be the fundamental weights:
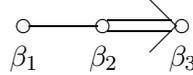
\begin{center}
\begin{picture}(100,80)
\put(10,50){\circle{5}}
\put(12,50){\line(1,0){26}}
\put(5,35){$\beta_1$}
\put(40,50){\circle{5}}
\put(35,35){$\beta_2$}
\put(42,52){\line(1,0){24}}
\put(68,50){\line(-1,1){10}}
\put(68,50){\line(-1,-1){10}}
\put(42,48){\line(1,0){24}}

\put(70,50){\circle{5}}
\put(65,35){$\beta_3$}
\end{picture}
\end{center}
Denote the nonzero weights of the representation $V(\pi_1)$ by $\pm\varepsilon_1, \pm\varepsilon_2, \pm\varepsilon_3$. One has:

$$\beta_1=\varepsilon_1-\varepsilon_2\quad
\beta_2=\varepsilon_2-\varepsilon_3\quad
\beta_3=\varepsilon_3,$$
$$\pi_1=\varepsilon_1\quad
\pi_2=\varepsilon_1+\varepsilon_2\quad
\pi_3=\frac{1}{2}(\varepsilon_1+\varepsilon_2+\varepsilon_3).$$

Let us number the positive roots of $B_3$ as follows:
\begin{center}
\begin{tabular}{lll}
 $\alpha_1=\varepsilon_1+\varepsilon_2$ &
 $\alpha_2=\varepsilon_1+\varepsilon_3$ &
 $\alpha_3=\varepsilon_2+\varepsilon_3$\\
 $\alpha_4=\varepsilon_1$ &
 $\alpha_5=\varepsilon_3$ &
 $\alpha_6=\varepsilon_2$ \\
 $\alpha_7=\varepsilon_1-\varepsilon_3$ &
 $\alpha_8=\varepsilon_2-\varepsilon_3$ &
 $\alpha_9=\varepsilon_1-\varepsilon_2$
\end{tabular}
\end{center}

We will use the following order on $\mathbb{Z}^9$: 
$\overline{p}=(p_1,\ldots,p_9)>\overline{q}=(q_1,\ldots,q_9)$ if $\sum p_i>\sum q_i$ and in case of a tie if $\overline{p}<\overline{q}$ lexicographically.

Restrict the representations of fundamental weights of $B_3$ to $G_2$:
\begin{center}
$
V_{B_3}(\pi_1)^*\mid_{G_2}= V_{G_2}(\omega_1)^*.
$\\

$
V_{B_3}(\pi_2)^*\mid_{G_2}=V_{G_2}(\omega_1)^*\oplus V_{G_2}(\omega_2)^*.
$\\
$
V_{B_3}(\pi_3)^*\mid_{G_2}=V_{G_2}(0)^*\oplus V_{G_2}(\omega_1)^*.
$
\end{center}
For every representation $V_{G_2}(\lambda)^*$ in $V_{B_3}(\pi_i)^*, i=1,2,3,$ we attach a signature $\sigma_{\pi_i;\lambda}$ of the lowest vector in $V_{G_2}(\lambda)^*$. One has 

\begin{flushleft}
\begin{tabular}{ll}
$\sigma_1:=(\pi_1;0,0,0,0,0,0,0,0,0)=\sigma_{\pi_1,\omega_1}$& $\sigma_4:=(\pi_3;0,0,0,0,0,0,0,0,0)=\sigma_{\pi_3,\omega_1}$\\
$\sigma_2:=(\pi_2;0,0,0,0,0,0,0,0,0)=\sigma_{\pi_2,\omega_2}$&$\sigma_5:=(\pi_3;0,0,1,0,0,0,0,0,0)=\sigma_{\pi_3,0}$\\
$\sigma_3:=(\pi_2;0,0,0,0,0,1,0,0,0)=\sigma_{\pi_2,\omega_1}$&
\end{tabular}

\end{flushleft}These signatures do not generate $\Sigma'$, because $d(\lambda)\neq V(\lambda)$ for $\lambda=\pi_1+\pi_2$ and $\lambda=\pi_1+\pi_3$. So we restrict these representations as well:
\begin{center}
$
V_{B_3}(\pi_1+\pi_2)^*\mid_{G_2}= V_{G_2}(\omega_1+\omega_2)^*\oplus V_{G_2}(2\omega_1)^*\oplus V_{G_2}(\omega_2)^*.
$\\

$
V_{B_3}(\pi_1+\pi_3)^*\mid_{G_2}=V_{G_2}(2\omega_1)^*\oplus V_{G_2}(\omega_1)^*\oplus V_{G_2}(\omega_2)^*.
$
\end{center}
This gives us two new signatures in $\Sigma'$:

\begin{flushleft}
\begin{tabular}{ll}
$\sigma_6:=(\pi_1+\pi_3;0,0,0,0,1,0,0,0,0)=\sigma_{\pi_1+\pi_3,\omega_2}$& \\
$\sigma_7:=(\pi_1+\pi_2;0,0,1,0,0,0,0,0,0)=\sigma_{\pi_1+\pi_2,\omega_2}$&
\end{tabular}
\end{flushleft}

It can be easily checked that $\sigma_i,i=1,\ldots,7,$ generate $\Sigma'$ with one relation $\sigma_4+\sigma_7=\sigma_1+\sigma_2+\sigma_5$. This solves the branching problem.

\subsection{The branching rule for $F_4\supset B_4$}

 We consider the standard regular embedding $B_4\subset F_4$.

Let $\beta_1, \beta_2, \beta_3, \beta_4$ be the  simple roots for $F_4$ and let $\pi_1, \pi_2, \pi_3, \pi_4$ be the fundamental weights:
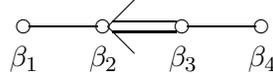
\begin{center}
\begin{picture}(100,80)
\put(10,50){\circle{5}}
\put(12,50){\line(1,0){26}}
\put(5,35){$\beta_1$}
\put(40,50){\circle{5}}
\put(35,35){$\beta_2$}
\put(44,52){\line(1,0){24}}
\put(42,50){\line(1,1){10}}
\put(42,50){\line(1,-1){10}}
\put(44,48){\line(1,0){24}}
\put(72,50){\line(1,0){26}}
\put(70,50){\circle{5}}
\put(65,35){$\beta_3$}
\put(100,50){\circle{5}}
\put(95,35){$\beta_4$}
\end{picture}
\end{center}
Denote the fundamental weights of $B_4$ by $\omega_1,\ldots,\omega_4$. Let the nonzero weights of the simplest ($\dim=9$) representation $V_{B_4}(\omega_1)$ be $\pm\varepsilon_1, \pm\varepsilon_2, \pm\varepsilon_3$. One has:

$$\beta_1=\frac{1}{2}(\varepsilon_1-\varepsilon_2-\varepsilon_3-\varepsilon_4),\quad
\beta_2=\varepsilon_4,\quad
\beta_3=\varepsilon_3-\varepsilon_4,\quad
\beta_4=\varepsilon_2-\varepsilon_3,$$
$$\pi_1=\varepsilon_1,\quad
\pi_2=\frac{1}{2}(3\varepsilon_1+\varepsilon_2+\varepsilon_3+\varepsilon_4),\quad
\pi_3=2\varepsilon_1+\varepsilon_2+\varepsilon_3,\quad
\pi_4=\varepsilon_1+\varepsilon_2.$$

Let us number the positive roots of $F_4$ as follows:
\begin{center}
$\alpha_1,\dots,\alpha_{16}=$ the roots of $B_4$ in any order,

\begin{tabular}{ll}
 $\alpha_{17}=\frac{1}{2}(\varepsilon_1+\varepsilon_2+\varepsilon_3+\varepsilon_4)$ &
 $\alpha_{18}=\frac{1}{2}(\varepsilon_1+\varepsilon_2+\varepsilon_3-\varepsilon_4)$\\
 $\alpha_{19}=\frac{1}{2}(\varepsilon_1+\varepsilon_2-\varepsilon_3-\varepsilon_4)$&
 $\alpha_{20}=\frac{1}{2}(\varepsilon_1+\varepsilon_2-\varepsilon_3+\varepsilon_4)$ \\
 $\alpha_{21}=\frac{1}{2}(\varepsilon_1-\varepsilon_2+\varepsilon_3-\varepsilon_4)$ &
 $\alpha_{22}=\frac{1}{2}(\varepsilon_1-\varepsilon_2+\varepsilon_3+\varepsilon_4)$ \\
 $\alpha_{23}=\frac{1}{2}(\varepsilon_1-\varepsilon_2-\varepsilon_3+\varepsilon_4)$ &
 $\alpha_{24}=\frac{1}{2}(\varepsilon_1-\varepsilon_2-\varepsilon_3-\varepsilon_4)$ \\

\end{tabular}
\end{center}

Now we introduce an order on $\mathbb{Z}^{24}$.
For the tuple $(p_1,\ldots, p_{24})$
 set
 $$
  q_{i}=\sum_{j=17}^{25-i}p_{j},
 $$
  Then $\sigma<\sigma'$ if $(q_{1},\ldots,q_{8})<(q_{1}',\ldots,q_{8}')$ in the lexicographic order, and in case of a tie we complete the order in any way. 
\begin{remark}
\emph{We are interested in a description of $\Sigma'$. It is easy to see that all signatures in $\Sigma'$ have zero exponents corresponding to the roots $\alpha_1,\ldots,\alpha_{16}$. Thus we do not care much about a completion of the order on $\mathbb{Z}^{24}$. }
\end{remark}

Restrict the representations of fundamental weights of $F_4$ to $B_4$:
\begin{center}
$
V_{F_4}(\pi_1)^*\mid_{B_4}= V_{B_4}(\omega_1)^*\oplus V_{B_4}(\omega_4)^*\oplus V_{B_4}(0)^*.
$\\

$
V_{F_4}(\pi_2)^*\mid_{B_4}= V_{B_4}(\omega_1+\omega_4)^*\oplus V_{B_4}(\omega_1)^*\oplus V_{B_4}(\omega_2)^*\oplus V_{B_4}(\omega_3)^*\oplus V_{B_4}(\omega_4)^*.
$\\
$
V_{F_4}(\pi_3)^*\mid_{B_4}= V_{B_4}(\omega_1+\omega_3)^*\oplus V_{B_4}(\omega_1+\omega_4)^*\oplus V_{B_4}(\omega_2+\omega_4)^*\oplus$\\
$\oplus V_{B_4}(\omega_3)^*\oplus V_{B_4}(\omega_2)^*.
$
\\
$
V_{F_4}(\pi_4)^*\mid_{B_4}= V_{B_4}(\omega_2)^*\oplus V_{B_4}(\omega_4)^*.
$
\end{center}

For every representation $V_{B_4}(\lambda)^*$ in $V_{F_4}(\pi_i)^*, i=1,2,3,4,$ we attach a signature $\sigma_{\pi_i;\lambda}$ of the lowest vector in $V_{B_4}(\lambda)^*$. We will omit the first 16 exponents of the signatures, because they are all zero (see the Remark above). One has: 

\begin{small}
\begin{tabular}{ll}
$\sigma_1:=(\pi_1;0,0,0,0,0,0,0,0)=\sigma_{\pi_1,\omega_1}$ &

$\sigma_9:=(\pi_3;0,0,0,0,0,0,0,0)=\sigma_{\pi_3,\omega_1+\omega_3}$ \\
$\sigma_2:=(\pi_1;1,0,0,0,0,0,0,0)=\sigma_{\pi_1,\omega_4}$ &

$\sigma_{10}:=(\pi_3;0,0,0,0,1,0,0,0)=\sigma_{\pi_3,\omega_1+\omega_4}$\\
$\sigma_3:=(\pi_1;1,0,1,0,0,0,0,0)=\sigma_{\pi_1,0}$ &

$\sigma_{11}:=(\pi_3;0,0,0,0,0,1,0,0)=\sigma_{\pi_3,\omega_2+\omega_4}$\\
$\sigma_4:=(\pi_2;0,0,0,0,0,0,0,0)=\sigma_{\pi_2,\omega_1+\omega_4}$ &

$\sigma_{12}:=(\pi_3;0,0,0,1,0,0,1,0)=\sigma_{\pi_3,\omega_3}$\\
$\sigma_5:=(\pi_2;0,0,1,0,0,0,0,0)=\sigma_{\pi_2,\omega_1}$ &

$\sigma_{13}:=(\pi_3;0,0,0,1,1,0,0,0)=\sigma_{\pi_3,\omega_2}$\\
$\sigma_6:=(\pi_2;0,0,0,1,0,0,0,0)=\sigma_{\pi_2,\omega_2}$ &

$\sigma_{14}:=(\pi_4;0,0,0,0,0,0,0,0)=\sigma_{\pi_4,\omega_2}$\\
$\sigma_7:=(\pi_2;0,1,0,0,0,0,0,0)=\sigma_{\pi_2,\omega_3}$ &

$\sigma_{15}:=(\pi_4;0,0,0,0,0,0,1,0)=\sigma_{\pi_4,\omega_4}$\\
$\sigma_8:=(\pi_2;0,1,0,0,1,0,0,0)=\sigma_{\pi_2,\omega_4}$ &\\
\end{tabular}
\end{small}
$\newline$
These signatures do not generate $\Sigma'$, because $d(\lambda)\neq V(\lambda)$ for $\lambda=\pi_1+\pi_3$,$\pi_1+\pi_4$, $\pi_2+\pi_4$, $\pi_3+\pi_4$. 
The restrictions to $B_4$ of these representations of $F_4$ give us five new signatures in $\Sigma'$:\\
\begin{small}
\begin{tabular}{l}
$
\sigma_{16}:=(\pi_1+\pi_4;1,0,0,0,0,0,0,1)=\sigma_{\pi_1+\pi_4,\omega_3}$\\
$\sigma_{17}:=(\pi_3+\pi_4;0,0,0,0,1,0,0,1)=\sigma_{\pi_3+\pi_4,\omega_1+\omega_3}$\\
$\sigma_{18}:=(\pi_2+\pi_4;0,0,0,0,0,0,0,1)=\sigma_{\pi_2+\pi_4,\omega_1+\omega_3}$\\
$\sigma_{19}:=(\pi_1+\pi_3;1,0,0,1,0,0,0,0)=\sigma_{\pi_1+\pi_3,\omega_2+\omega_3}$\\
$\sigma_{20}:=(\pi_2+\pi_4;0,1,0,0,1,0,0,1)=\sigma_{\pi_2+\pi_4,\omega_3}$
\end{tabular}
\end{small}\\
It can be easily checked by using a computer that $\sigma_i,i=1,\ldots,20,$ generate $\Sigma'$ with 28 relations given below. This solves the branching problem. These relations form the reduced Groebner basis of the ideal of all relations with respect to the lexicographic order ($\sigma_i>\sigma_j$ if $i<j$).  The relations are:\\
\begin{small}
\begin{tabular}{ll}
$1.\,\sigma_{10}+\sigma_{12}+\sigma_{14}+\sigma_{16}=\sigma_{15}+\sigma_{17}+\sigma_{19}$
& 

$20.\, \sigma_3+\sigma_6+\sigma_{9}=\sigma_{5}+\sigma_{19}$\\
$2.\, \sigma_{9}+\sigma_{13}+\sigma_{16}=\sigma_{17}+\sigma_{19}$
& 

$21.\, \sigma_2+\sigma_{20}=\sigma_{8}+\sigma_{16}$\\
$3.\, \sigma_{9}+\sigma_{13}+\sigma_{15}=\sigma_{10}+\sigma_{12}+\sigma_{14}$
& 

$22.\, \sigma_2+\sigma_{18}=\sigma_{4}+\sigma_{16}$\\
$4.\, \sigma_{8}+\sigma_{17}=\sigma_{10}+\sigma_{20}$&

$23.\, \sigma_2+\sigma_{17}=\sigma_{10}+\sigma_{16}$\\
$5.\, \sigma_{8}+\sigma_{12}+\sigma_{14}+\sigma_{16}=\sigma_{15}+\sigma_{19}+\sigma_{20}$&

$24.\, \sigma_2+\sigma_{12}+\sigma_{14}=\sigma_{15}+\sigma_{19}$\\

$6.\, \sigma_{7}+\sigma_{17}=\sigma_{9}+\sigma_{20}$&
$25.\, \sigma_2+\sigma_{9}+\sigma_{13}=\sigma_{10}+\sigma_{19}$\\
$7.\, \sigma_{7}+\sigma_{13}+\sigma_{16}=\sigma_{19}+\sigma_{20}$&

$26.\,\sigma_2+\sigma_{7}+\sigma_{13}=\sigma_{8}+\sigma_{19}$\\
$8.\, \sigma_{7}+\sigma_{13}+\sigma_{15}=\sigma_8+\sigma_{12}+\sigma_{14}$&

$27.\, \sigma_2+\sigma_{6}+\sigma_{9}=\sigma_{4}+\sigma_{19}$\\
$9.\, \sigma_{7}+\sigma_{10}=\sigma_{8}+\sigma_{9}$&

$28.\, \sigma_2+\sigma_{5}=\sigma_{3}+\sigma_{4}$\\
$10.\, \sigma_{6}+\sigma_{17}=\sigma_{13}+\sigma_{18}$&\\
$11.\, \sigma_{6}+\sigma_{10}+\sigma_{20}=\sigma_8+\sigma_{13}+\sigma_{18}$&\\
$12.\, \sigma_{6}+\sigma_{9}+\sigma_{20}=\sigma_7+\sigma_{13}+\sigma_{18}$&\\
$13.\, \sigma_{6}+\sigma_9+\sigma_{16}=\sigma_{18}+\sigma_{19}$&\\
$14.\, \sigma_{4}+\sigma_{20}=\sigma_{8}+\sigma_{18}$&\\
$15.\, \sigma_4+\sigma_{17}=\sigma_{10}+\sigma_{18}$&\\
$16.\, \sigma_4+\sigma_{13}=\sigma_{6}+\sigma_{10}$&\\
$17.\, \sigma_4+\sigma_{12}+\sigma_{14}=\sigma_{6}+\sigma_{9}+\sigma_{15}$&\\
$18.\, \sigma_3+\sigma_{18}=\sigma_{5}+\sigma_{16}$&\\
$19.\, \sigma_{3}+\sigma_{6}+\sigma_{10}+\sigma_{12}+\sigma_{14}=\sigma_{5}+\sigma_{13}+\sigma_{15}+\sigma_{19}$&\\
\end{tabular}
\end{small}

\end{document}